\documentclass{amsart}

\usepackage[all]{xy}

\theoremstyle{plain}
\newtheorem{theorem}{Theorem}[section]
\newtheorem{proposition}[theorem]{Proposition}
\newtheorem{lemma}[theorem]{Lemma}
\newtheorem{corollary}[theorem]{Corollary}
\newtheorem{conjecture}[theorem]{Conjecture}

\newtheorem{question}[theorem]{Question}

\theoremstyle{definition}
\newtheorem{definition}[theorem]{Definition}
\newtheorem{example}[theorem]{Example}

\theoremstyle{remark}
\newtheorem{remark}[theorem]{Remark}
\newtheorem{remarks}[theorem]{Remarks}
\newtheorem{discussion}[theorem]{Discussion}

\newcommand{\thmref}[1]{Theorem~\ref{#1}}
\newcommand{\propref}[1]{Proposition~\ref{#1}}
\newcommand{\lemref}[1]{Lemma~\ref{#1}}
\newcommand{\corref}[1]{Corollary~\ref{#1}}
\newcommand{\conjref}[1]{Conjecture~\ref{#1}}

\newcommand{\exref}[1]{Example~\ref{#1}}
\newcommand{\queref}[1]{Question~\ref{#1}}
\newcommand{\remref}[1]{Remark~\ref{#1}}
\newcommand{\remsref}[1]{Remarks~\ref{#1}}

\def\Lam{\Lambda}
\def\set#1{\{#1 \}}
\def\a{\alpha}
\def\odd{\mathrm{odd}}
\def\rank{\mathrm{dim}}
\def\sph{{\mathbf S}}
\def\id{\operatorname{ id}}

\begin{document}

\title{Free Torus Actions and Two-Stage Spaces}

\author{Barry Jessup}

\address{Department of Mathematics,
          University of Ottawa,
          Ottawa
          Canada K1N 6N5}

\email{Bjessup@sciences.uottawa.ca}
\thanks{This research was supported in part by NSERC}
\author{Gregory Lupton}

\address{Department of Mathematics,
          Cleveland State University,
          Cleveland OH 44115
          U.S.A.}

\email{Lupton@math.csuohio.edu}

\date{\today}



\begin{abstract}
We prove  the toral rank conjecture of Halperin in some new cases.
Our results apply to certain  elliptic spaces that have a
two-stage Sullivan minimal model, and are obtained by combining
new lower bounds for the dimension of the cohomology  and new
upper bounds for the
  toral rank.  The paper
concludes with  examples and suggestions for future work.
\end{abstract}

\maketitle

\section{Introduction}

A well-known conjecture due to Halperin concerns torus actions on
a space (see \cite[Prob.1.4]{Halperin}). If $X$ is a space on
which an $n$-dimensional torus acts,  we say the action is
\emph{almost-free} if each isotropy subgroup is a finite group.
The largest integer $n \geq 1$ for which $X$ admits an almost-free
$n$-torus action is called the \emph{toral rank} of $X$, and is
denoted  $\mathrm{rk}\,(X)$.  If $X$ does not admit an almost-free
circle action, then   $\mathrm{rk}\,(X) = 0$. Halperin's
conjecture gives an upper bound for the toral rank of $X$ in terms
of its cohomology, as follows:

\begin{conjecture}[The toral rank conjecture] \label{conj:TRC}
If $X$ is simply connected, then
$$\mathrm{dim}\,H(X; \mathbb{Q}) \geq 2^{\mathrm{rk}\,(X)}.$$
\end{conjecture}

We shall henceforth assume that all our spaces are 1-connected,
finite cell complexes.  There are some technical conditions on the
topology of the space $X$ in Halperin's original formulation, but
as these are satisfied by finite cell complexes, we will not
mention them explicitly here.

The main tool we shall use is the Sullivan \emph{minimal model},
and a   basic reference is \cite{F-H-T00}. For our purposes, we
note that to any 1-connected space $X$ there corresponds, in a
contravariant way, a commutative differential graded algebra
(unique up to isomorphism) $(\Lambda W,d)$, called the minimal
model of $X$, which algebraically models the rational homotopy
type of the space. By $\Lambda W$ we mean the free graded
commutative algebra generated by the graded vector space $W$. The
differential $d$ of any element of $W$ is a polynomial in $\Lambda
W$ with no linear term, which in particular means that there is  a
homogeneous basis  $\{w_i\}_{i\ge 1}$ of $W$   for which
$dw_i\in\Lam W_{<i}$, where  $W_{<i}$ denotes the subspace of $W$
generated by $\{w_j\}_{j<i}$. (We will also occasionally consider
models $(\Lam W,d)$ which satisfy this latter nilpotence condition
but where the polynomial $dw$ may have a linear term.) This
contravariant correspondence yields an equivalence between the
homotopy category of 1-connected  rational spaces of finite type
and that of 1-connected  rational commutative differential graded
algebras of finite type. In particular, if $(\Lambda W,d)$ is the
minimal model of $X$, then $H(\Lambda W,d) \cong H(X; \mathbb{Q})$
as graded  algebras, and $W  \cong  \pi (X)\otimes \mathbb{Q}$ as
graded vector spaces. Moreover, we remark that  a series of
results of Halperin and others implies that if $H(\Lambda W,d)$
and $W$ are finite-dimensional, then $H(\Lambda W,d)$ must satisfy
Poincar{\'e} duality, and $(\Lambda W,d)$ is the minimal model of
some closed smooth manifold. In all cases of interest to us here,
the space of generators $W$ will be finite-dimensional.  If
$\{w_1, \ldots, w_m\}$ is a homogeneous basis for $W$, then we
will write $W = \langle w_1, \ldots, w_m \rangle$, and we also
denote $\Lambda W$ by $\Lambda(w_1, \ldots, w_{m})$ in this case,
often omitting explicit reference to the differential. We note
that models with $W$  and $H(\Lambda W,d)$ both finite-dimensional
are called {\it elliptic}, and   a space $X$ with an   elliptic
minimal model is called an elliptic space. Topologically, this
means that both $\pi(X)\otimes \mathbb{Q}$ and $H(X; \mathbb{Q})$
are finite dimensional.

\conjref{conj:TRC} is already known to hold for homogeneous spaces
$G/H$, for $G$ connected, and $H$ closed and connected
\cite[Prop.1.5]{Halperin}. Such spaces
  have
\emph{two-stage}  minimal models \cite[Prop.15.16]{F-H-T00}, where
a minimal model $(\Lambda W, d)$ is said to be  two-stage if $W$
decomposes as $W \cong U \oplus V$ with $d\,U  = 0$ and $d\, V
\subseteq \Lambda U$. By the remark in the previous paragraph, it
is easy to see that there are many other examples of spaces,
indeed smooth manifolds, with two-stage models, and it is these
spaces that we shall study (see for example \corref{cor:TRC
pure}).

  We note for later reference that there may
be several ways to display a minimal model as a two-stage model.
In particular, a generator in $V$ that is a cocycle could just as
well be included in $U$.  We will generally be interested in
choosing a two-stage decomposition
   in which $V$ is as large as possible.

By a K-S extension (Koszul-Sullivan), or simply an extension, we
mean a sequence of the form
\begin{displaymath}
(\Lambda W_1, d_1) \to (\Lambda W_1 \otimes \Lambda W_2, D) \to
(\Lambda W_2, d_2)
\end{displaymath}
with $(\Lambda W_1, d_1)$ a minimal model, in which $D$ restricts
to $d_1$ on $\Lambda W_1 \otimes 1$, and for which there is an
ordered basis of $W_2 = \langle w_1, w_2, \ldots \rangle$ with $D
w_i  \in \Lambda W_1 \otimes \Lambda(w_1, \ldots, w_{i-1})$ for
each $i$. A K-S extension is the minimal-model analogue of a Serre
fibration (cf.~\cite[Sec.15(a)]{F-H-T00}), and for this reason,
$\Lambda W_1, \Lambda W_1 \otimes \Lambda W_2 $ and $\Lambda W_2$
are known respectively as the base, the total space and the fibre
of the extension.

The connection between minimal models and the toral rank is
originally due to Allday and Halperin \cite{Allday-Halperin}
(cf.~\cite[Prop.4.2]{Halperin}): If $X$ has minimal model
$(\Lambda W, d)$ and  admits an $n$-dimensional torus action, then
there is an extension
\begin{equation}\label{eq:Borel fibration}
\Lambda A_n \to  \Lambda A_n \otimes \Lambda W  \to  \Lambda W
\end{equation}
in which $A_n = \langle a_1, \dots, a_n \rangle$ with each $a_i$
of degree $2$.  If the action is almost free, then
$\mathrm{dim}\,H(\Lambda A_n \otimes \Lambda W, D)$ is
\emph{finite-dimensional} \cite[Prop.4.2]{Halperin}. In principle,
this result allows an upper bound for $\mathrm{rk}\,(X)$ to be
obtained by a direct analysis of the minimal model of $X$, and
this direct approach has been carried out to great effect in some
situations. In general, however, the computational problems
involved here appear to be quite substantial.

However, since we are really only considering the rational
homotopy type of $X$, we are led to the following variation of the
toral rank: the \emph{rational toral rank} $\mathrm{rk}_0(X)$ of
$X$ is defined by $\mathrm{rk}_0(X) = \text{max}\{\mathrm{rk}(Y)
\mid Y \simeq_{\mathbb{Q}} X \}$.  Clearly, we have
$\mathrm{rk}_0(X) \geq \mathrm{rk}(X)$, hence an upper bound on
$\mathrm{rk}_0(X)$ will serve as one on $\mathrm{rk}(X)$. See
\cite[Prop.4.2]{Halperin} for the precise relationship between
these two numbers.    The characterization of $\mathrm{rk}_0(X)$
in terms of a minimal model of $X$ is as follows. If  $(\Lambda W,
d)$ is a minimal model of $X$, then $\mathrm{rk}_0(X)$  is the
largest $n$ (if such exists) for which there is a K-S-extension of
the form (\ref{eq:Borel fibration}), for which
$\mathrm{dim}\,H(\Lambda A_n \otimes \Lambda W, D)$ is
finite-dimensional.  We will also denote $\mathrm{rk}_0(X)$ by
$\mathrm{rk}_0(\Lambda W)$.

Now we briefly summarize the contents of the paper. In Section 2,
\thmref{thm:general 2-stage}  gives a lower bound for the
dimension of the cohomology of any space with an elliptic
two-stage minimal model. This result is proved using standard
tools familiar from algebraic topology, namely   the \emph{Wang
sequence} and the Serre spectral sequence. Also in Section 2 we
give two results that, under certain additional hypotheses, give
upper bounds on the rational toral rank---hence on the toral
rank---of a space with two-stage minimal model (\thmref{thm:two
stage even and odd} and \thmref{thm:stable and separated}).
Combining these results allows us to establish \conjref{conj:TRC}
in some new cases (\corref{cor:two stage even and odd}).
  In addition, our results can be used to provide new and
essentially elementary proofs of the conjecture in some
previously-known cases, which we   discuss  after
\thmref{thm:general 2-stage}. In Section 3, we give a number of
examples and suggest several directions for future investigation.

\subsection*{Acknowledgements} It is our pleasure to thank
Yves F{\'e}lix and Steve Halperin for their input. They have each
generously shared with us their expertise in concocting examples
of the sort germane to this paper. In particular,
\exref{ex:Steve's example} and the discussion that follows it are
entirely due to Steve Halperin.

\section{Cohomology and Toral Rank of Two-Stage Minimal Models}%
\label{sec:2-stage}

In this section we give several results about the cohomology and
  the toral rank of two-stage minimal models.

\subsection{Cohomology of two-stage minimal models}
We begin by establishing a lower bound for the dimension of the
cohomology of an elliptic space with two-stage minimal model.  Our
bound is of a form similar to that featured in the toral rank
conjecture.

  If $\Lambda(U, V)$ is an
elliptic two-stage minimal model, then the Mapping theorem
\cite[P.375]{F-H-T00} implies that $V$ has generators of odd
degree only. On the other hand, $U$ may have generators of odd or
even degree. For the remainder of the paper, unless otherwise
specified, we assume that all minimal models are two-stage,
elliptic.  For details about elliptic spaces and their minimal
models, see \cite[Sec.32]{F-H-T00}.

First, we establish our lower bound in a special case.  The main
result is then established by reducing to this special case.

\begin{proposition}\label{prop:odd-free}
Let $\Lambda(U, V)$ be a two-stage, elliptic minimal model with
odd degree generators only and suppose that $d \colon V \to
\Lambda^2U$ is an isomorphism. Then
$\mathrm{dim}\,H\big(\Lambda(U, V)\big) \geq 2^{\mathrm{dim}\,V}$.
\end{proposition}

\begin{proof}  Suppose that $\mathrm{dim}\,U = n$.  Since
$d \colon V \to \Lambda^2U$ an isomorphism, and since $\Lambda U$
is an exterior algebra, we have $\mathrm{dim}\,V = \binom{n}{2}$.
Thus we must show that $\mathrm{dim}\,H\big(\Lambda(U, V)\big)
\geq 2^{\binom{n}{2}}$. We will proceed by induction on $n$.
First, we give an explicit description of the model $\Lambda(U,
V)$.

Write $\Lambda(U, V)$ as ${\mathcal M}_n = (\Lambda(U_n, V_n),
d_n)$, with $U_n = \langle u_1, \ldots, u_n \rangle$ and $V_n =
\langle \{v_{i,j} \}_{1 \leq i < j \leq n} \rangle$.  The
differential is $d_{n}(U_{n}) = 0$ and $d_{n}(v_{i, j}) = u_iu_j$
for each $0 \leq i < j \leq n$. For $n = 1$, we have
$\mathrm{dim}\, U_1  = 1$ and $V_1 = 0$, and   the proposition is
trivial. For $n = 2$, we have ${\mathcal M}_2 = \Lambda(u_1, u_2,
v_{1, 2})$. It is easily checked that $\mathrm{dim}\,H({\mathcal
M}_2) = 6 \geq 2^1$. This starts the induction.

Now suppose inductively that $\mathrm{dim}\,H({\mathcal M}_n) \geq
2^{\binom{n}{2}}$ for some $n \geq 2$.  We adjust the notation
somewhat and write ${\mathcal M}_{n+1}$ as follows.  Write
$U_{n+1} = \langle u_0, u_1, \ldots, u_n \rangle = \langle
u_0\rangle \oplus U_n$.  Also, set $V_0 = \langle v_{0, 1},
\ldots, v_{0, n} \rangle$, so that $V_{n+1} = V_0 \oplus V_n$. The
differential $d_{n+1}$ of ${\mathcal M}_{n+1}$ extends that of
${\mathcal M}_n$. Thus, we have $d_{n+1}(U_{n+1}) = 0$,
$d_{n+1}(v_{i, j}) = d_{n}(v_{i, j}) = u_iu_j$ for $1 \leq i < j
\leq n$ and $d_{n+1}(v_{0, j}) = u_0u_j$ for $1 \leq j \leq n$.

Further, if $\overline{\mathcal{M}}_{n+1}$ denotes $ {\mathcal
M}_{n}\otimes \Lambda V_0$,  we have the following extension:
\begin{equation}\label{eq:Mn extension}
\Lambda  u_0  \to ({\mathcal M}_{n+1}, d_{n+1}) \to
(\overline{\mathcal{M}}_{n+1}, \overline{d}_{n+1}).
\end{equation}

Consider the short exact sequence of differential graded vector
spaces
\begin{equation}\label{eq:exact vector space}
\xymatrix{ 0 \ar[r] & \overline{\mathcal{M}}_{n+1} \ar[r]^-{j} &
\mathcal{M}_{n+1} \ar[r]^-{p} &
  \overline{\mathcal{M}}_{n+1} \ar[r] & 0, \\}
\end{equation}
where $p$ denotes the projection and $j$ the map defined by
$j(\chi) = u_0\,\chi$ for $\chi \in \Lambda W$. The ensuing long
exact sequence in cohomology has a connecting homomorphism that
may be described as follows: On $\overline{\mathcal{M}}_{n+1}$,
define $\theta^*:\overline{\mathcal{M}}_{n+1}\to
\overline{\mathcal{M}}_{n+1}$ by $u_0\theta=d_{n+1}-\bar d_{n+1}
$. Then, one can easily check that $\theta$ is  actually a
derivation of the algebra $\overline{\mathcal{M}}_{n+1}$ of degree
$1- |u_0|$ which commutes with $\bar d_{n+1}$, and so induces a
map $\theta^* : H(\overline{\mathcal{M}}_{n+1}) \to H^{*-
|u_0|+1}(\overline{\mathcal{M}}_{n+1})$ which is the connecting
homomorphism of (\ref{eq:exact vector space}). We will call the
derivation $\theta^*$ the \emph{Wang derivation} for (\ref{eq:Mn
extension}), for reasons indicated in \remref{wang derivation}
below.

This Wang derivation can be explicitly computed as follows.
Because  $d_{n+1}(V_0) \subseteq u_0\cdot U_n$, we see  that
$H(\overline{\mathcal{M}}_{n+1},\overline{d}_{n+1})=H({\mathcal
M}_{n})\otimes \Lambda V_0$. Now, for $\chi \in {\mathcal M}_{n}
\otimes \Lambda V_0 $,  we find that $\theta({\mathcal M}_{n}) =
0$, directly from the definition. Since $d_{n+1}(v_{0, i}) =
u_0u_i$, we have $\theta(v_{0, i}) = u_i$ for each $i = 1, \ldots,
n$.  On passing to cohomology, therefore, we find that
$\theta^*\big(H({\mathcal M}_{n})\big) = 0$ and $\theta^*(v_{0,
i}) = u_i \in H({\mathcal M}_{n})$ for $i = 1, \ldots, n$.

As usual, we can condense the long exact sequence coming from
(\ref{eq:exact vector space}) into
\begin{displaymath}
0 \to \text{coker}\,\theta^* \to H({\mathcal M}_{n+1}) \to
\text{ker}\,\theta^* \to 0,
\end{displaymath}
so that to complete the induction, it simply remains to show that
$$\dim \text{ker}\,\theta^* +\dim \text{coker}\,\theta^*=2 \dim \ker \theta^*$$
is at least $2^{\binom{n+1}{2}}$.  Our inductive assumption is
that $\mathrm{dim}\,H({\mathcal M}_n) \geq 2^{\binom{n}{2}}$.
Therefore $\mathrm{dim}\,\big(H({\mathcal M}_n)\cdot \Lambda
V_0\big) \geq 2^{\binom{n}{2}}\times 2^n$. Since $\binom{n}{2} + n
= \binom{n+1}{2}$, our inductive assumption implies that
$\mathrm{dim}\,\big(H({\mathcal M}_n)\otimes \Lambda V_0\big) \geq
2^{\binom{n+1}{2}}$.

We claim that $(\theta^*)^2 = 0$, after which elementary linear
algebra completes the proof. The key observation here is the
following:  Denote by $[u_i] \in H({\mathcal M}_{n})$ the class
represented by the cocycle $u_i$.  Then the ideal of $H({\mathcal
M}_{n})$ generated by all two-fold products $\{ [u_i] [u_j] \}_{1
\leq i < j \leq n}$ is zero.  This follows since each product $u_i
u_j$ is a boundary, namely  $d_{n}(v_{i, j})$. Recall that
$\theta^*\big(H({\mathcal M}_n)\big) = 0$, whilst $\theta^*(V_0) =
U_n \subseteq H({\mathcal M}_n)$.  Thus,
$(\theta^*)^2\big(H({\mathcal M}_n)\otimes \Lambda V_0\big)$ is
contained in the ideal generated by the two-fold products $\{
[u_i] [u_j] \}_{1 \leq i < j \leq n}$, which we observed is zero.
We conclude that $(\theta^*)^2 = 0$.

But $(\theta^*)^2 = 0$ implies $ 2\dim \text{ker}\,\theta^* \geq
\mathrm{dim}\,\big(H({\mathcal M}_n)\otimes \Lambda V_0\big)$.
Therefore, we have
\begin{displaymath}
\mathrm{dim}\,H({\mathcal M}_{n+1})=  2\dim
\text{ker}\,\theta^*\geq \mathrm{dim}\big(H({\mathcal
M}_{n})\otimes\Lambda V_0\big) \geq 2^{\binom{n+1}{2}}.
\end{displaymath}
This completes the induction.
\end{proof}

\begin{remark}\label{wang derivation}
Topologically, the extension sequence (\ref{eq:Mn extension})
corresponds to a fibration with base an odd sphere $S^{2r+1}$, and
the long exact sequence in cohomology induced from (\ref{eq:exact
vector space})  corresponds to the Wang sequence of the fibration,
in the usual sense of a fibration with base a sphere. The
derivation $\theta^*$ corresponds to the Wang derivation of the
fibration, again in the usual sense (cf. \cite[p.319]{Whitehead}).
\end{remark}

As promised, the main result is now obtained by reducing the
general two-stage case to that of \propref{prop:odd-free}.

\begin{theorem}\label{thm:general 2-stage}
Suppose $\Lambda(U, V)$ is any two-stage, elliptic minimal model.
Then
\begin{displaymath}
\mathrm{dim}\,H\big(\Lambda(U, V)\big) \geq 2^{\mathrm{dim}\,V -
\mathrm{dim}\,U^{\mathrm{even}}}.
\end{displaymath}
\end{theorem}

\begin{proof}
Given $\Lambda(U, V)$, form the following extension sequence:
\begin{equation}\label{eq:extn1}
\Lambda(U, V) \to \Lambda(U, V)\otimes \Lambda(\overline{W}, W), D
\to \Lambda (\overline{W}, W), \overline{D} = 0,
\end{equation}
in which $D \colon \overline{W} \to U^{\mathrm{even}}$ and $d
\colon W \to \Lambda^2 (U^{\mathrm{odd}})$ are vector space
isomorphisms. This extension sequence has an associated  Serre
spectral sequence, obtained by filtering the  total space  by
degree in $\Lambda(U, V)$. This spectral sequence has $E_2$-term
isomorphic to $H(\Lambda(U, V), d)\otimes H(\Lambda (\overline{W},
W), \overline{D}=0) \cong H(\Lambda(U, V), d)\otimes \Lambda
(\overline{W}, W)$ and it converges to $H(\Lambda(U, V)\otimes
\Lambda(\overline{W}, W), D)$ (cf.~\cite[Sec.18]{F-H-T00},
especially Example 2 and Exercise 2 of that section). Since the
$E_2$-term of this spectral sequence must be at least as large as
the term that it converges to, we have the following inequality:
\begin{equation}\label{eq:1stineq}
\mathrm{dim}\big(H(\Lambda(U, V), d)\big)\cdot
2^{\mathrm{dim}\,\overline{W}} \cdot 2^{\mathrm{dim}\,W} \geq
\mathrm{dim}\big(H(\Lambda(U, V)\otimes \Lambda(\overline{W}, W),
D)\big).
\end{equation}

We now claim that $(\Lambda(U, V)\otimes \Lambda(\overline{W}, W),
D)$ and $(\Lambda(U^{\mathrm{odd}}, W)\otimes \Lambda(V),
D\otimes1)$ are quasi-isomorphic.  To see this, argue as follows.
First, display the middle term of (\ref{eq:extn1}) as the middle
term of the following extension sequence:
\begin{equation}\label{eq:extn2}
(\Lambda(U^{\mathrm{even}}, \overline{W}), D) \to (\Lambda(U,
V)\otimes \Lambda(\overline{W}, W), D) \to
(\Lambda(U^{\mathrm{odd}}, W, V), \overline{D}),
\end{equation}
in which $(\Lambda(U^{\mathrm{even}}, \overline{W}), D)$ is an
acyclic model.  It follows that the projection $(\Lambda(U,
V)\otimes \Lambda(\overline{W}, W), D) \to
(\Lambda(U^{\mathrm{odd}}, W, V), \overline{D})$ is a
quasi-isomorphism.  Next, define a map of models $\phi \colon
(\Lambda(U^{\mathrm{odd}}, W)\otimes \Lambda(V), D\otimes1) \to
(\Lambda(U^{\mathrm{odd}}, W, V), \overline{D})$ as follows.  Set
$\phi$ to be the identity on $\Lambda(U^{\mathrm{odd}}, W)$. For
each generator $v_i \in V$, we have $\overline{D}(v_i) \in
\Lambda^{\geq 2} U^{\mathrm{odd}} \subseteq
\overline{D}(W)\cdot\Lambda(U^{\mathrm{odd}})$, since
$\overline{D}(W) \cong \Lambda^2(U^{\mathrm{odd}})$.  So choose an
element $x_i \in W\cdot\Lambda(U^{\mathrm{odd}})$, such that
$\overline{D}(x_i) = \overline{D}(v_i)$, for each $i$.  Finally,
set $\phi(v_i) = v_i - x_i$ and extend to an algebra map.  Notice
that $\phi$ makes the following  diagram of extension sequences
commute:
\begin{displaymath}
\xymatrix{(\Lambda(U^{\mathrm{odd}}, W), \overline{D}) \ar[r]
\ar[d]_{1}& (\Lambda(U^{\mathrm{odd}}, W), \overline{D}) \otimes
\Lambda V, D\otimes 1 \ar[r] \ar[d]^{\phi} & \Lambda V,
\overline{D\otimes 1} =
0 \ar[d]^{1} \\
(\Lambda(U^{\mathrm{odd}}, W), \overline{D}) \ar[r]&
(\Lambda(U^{\mathrm{odd}}, W, V), \overline{D}) \ar[r] & \Lambda
V, \overline{D} = 0 }
\end{displaymath}
It follows that $\phi$ is a quasi-isomorphism.

Returning to the inequality (\ref{eq:1stineq}) displayed above, we
now have
\begin{displaymath}
\mathrm{dim}\,H\big(\Lambda(U, V)\big)\cdot
2^{\mathrm{dim}\,\overline{W}} \cdot 2^{\mathrm{dim}\,W} \geq
\mathrm{dim}\,H\big(\Lambda(U^{\mathrm{odd}}, W)\big) \cdot
2^{\mathrm{dim}\,V}.
\end{displaymath}
But $\mathrm{dim}\,H\big(\Lambda(U^{\mathrm{odd}}, W)\big) \geq
2^{\mathrm{dim}\,W}$, by \propref{prop:odd-free}. Since
$\mathrm{dim}\,\overline{W} = \mathrm{dim}\,U^{\mathrm{even}}$,
the result follows.
\end{proof}

We can establish \conjref{conj:TRC} in several cases using
\thmref{thm:general 2-stage}. We begin with a discussion of two
previously-known cases of \conjref{conj:TRC}. We include these
cases here to illustrate that \thmref{thm:general 2-stage} can be
used to unify them with other results. Also, we feel it worthwhile
to include a self-contained proof in the second case that uses
techniques from well within rational homotopy theory.

\subsection{The pure case} For any
elliptic minimal model $\Lambda W$, the \emph{homotopy Euler
characteristic} is
$\chi_{\pi}:=\mathrm{dim}(W^{\mathrm{even}})-\mathrm{dim}(W^{\mathrm{odd}})$.
For any  elliptic space, a result due to Allday and Halperin
\cite[Th.1]{Allday-Halperin} asserts that the rational toral rank
is less than or equal to the negative of the homotopy Euler
characteristic, thus $\mathrm{rk}_0(X) \leq -\chi_{\pi}$.

Now consider the so-called \emph{pure} case, in which $\Lambda(U,
V)$ is a two-stage, elliptic minimal model with $U =
U^{\mathrm{even}}$ and $V = V^{\mathrm{odd}}$.  In this case,
\thmref{thm:general 2-stage} specializes to:
$\mathrm{dim}\,H\big(\Lambda(U, V)\big) \geq 2^{-\chi_{\pi}}$.
Combining these observations, we retrieve the following result of
Halperin:

\begin{corollary}[{\cite[Prop.1.5]{Halperin}}]\label{cor:TRC pure}
Let $X$ be an elliptic space with a pure minimal model. Then
\conjref{conj:TRC} holds for $X$, i.e.,
$\mathrm{dim}\,H(X;\mathbb{Q}) \geq 2^{\mathrm{rk}_0(X)}$.
\end{corollary}

We remark that this includes the case in which $X$ is a
homogeneous space. Indeed, our proof of \thmref{thm:general
2-stage} incorporates the argument used by Halperin to obtain
\conjref{conj:TRC} in the homogeneous space case.

\subsection{The case of odd generators with
quadratic differential (coformal  spaces)}    We say that a
minimal model is coformal if it has  a \emph{quadratic
differential}, that is, if $d \colon W \to \Lambda^2W$. Notice
that this really means quadratic with respect to a particular
choice of generators for $\Lambda W$. A space is coformal if it
has a coformal minimal model, and in this case its rational
homotopy type is completely determined by its rational homotopy
Lie algebra, the bracket of which is dual to $d$ \cite[\S
21]{F-H-T00}.

The following result is due to Allday and Puppe:

\begin{theorem}[{\cite[Prop.3.1,
Cor.3.5]{Allday-Puppe85a}}]\label{thm:two stage quadratic} Let $X$
be a simply-connected, coformal elliptic space that has a
two-stage minimal model $\Lambda(U, V)$ with odd-degree generators
only. Assume that the two-stage decomposition displays $V$ with
maximal dimension.  Then, $\mathrm{rk}_0(X) = \mathrm{dim}\,V $
and \conjref{conj:TRC} holds for $X$.
\end{theorem}

\begin{proof}
The equality $\mathrm{rk}_0(X) = \mathrm{dim}\, V $ is given in
\cite[Cor.3.5]{Allday-Puppe85a}. (Actually, the rational toral
rank is identified there with the dimension of the centre of the
homotopy Lie algebra of $X$. However, standard ideas from rational
homotopy theory can be used to identify $\mathrm{dim} V $ and the
dimension of the centre of the homotopy Lie algebra, in this
particular case.) Since the minimal model has odd-degree
generators only, \thmref{thm:general 2-stage} specializes to the
inequality $\mathrm{dim}\,H\big(\Lambda(U, V)\big) \geq
2^{\mathrm{dim}\,V}$.
\end{proof}

This paper contains  a self-contained proof of \thmref{thm:two
stage quadratic} that uses only standard tools from rational
homotopy theory.  For this, we need to supply our own argument for
the equality $\mathrm{rk}_0(X) = \mathrm{dim}\,V$, to substitute
for \cite[Cor.3.5]{Allday-Puppe85a}. This is obtained by combining
\lemref{lem:upper bound} with \thmref{thm:two stage even and odd}
below, both of which apply to more general settings.
(cf.~\remsref{rems:special cases}). When we prove \thmref{thm:two
stage even and odd}, we will be more specific as to what the
phrase `displays $V$ with maximal dimension' entails.

We remark that Allday-Puppe conclude their result by appealing to
a result of Deninger-Singhof \cite{Den-Sing}.  Indeed, the result
from \cite{Den-Sing} used to complete their argument is precisely
a special case of \thmref{thm:general 2-stage}, in which the
differential is assumed to be quadratic and all generators are
assumed to be of odd degree. It is therefore natural to ask
whether $\mathrm{rk}_0\big(\Lambda(U, V)\big) = \mathrm{dim}\,V$
more generally, at least in the case in which the two-stage,
elliptic minimal model has odd-degree generators only. This is not
true without some further hypotheses, as we illustrate in
\exref{ex:two-stage odd}. We return to this point, and related
questions, in the last section of the paper.

\subsection{New cases} In the remainder of this section,
we prove new results that give further situations to which we can
apply \thmref{thm:general 2-stage}. In phrasing our results, we
must be careful about a technical point, namely, the choice of a
two-stage decomposition that displays $V$ with maximal dimension.
We illustrate this  with an example.

\begin{example}\label{ex:non-maximal V}
Let $\Lambda(U, V; d)$ be a two-stage minimal model with bases $U
= \langle u_1, u_2, u_3 \rangle$ and $V = \langle v_1 \rangle$,
all generators being of degree $1$. If we set $d v_1  = u_1 u_3 -
u_1 u_2 + u_2u_3$, then each basis element of $U$ certainly occurs
in a non-trivial differential. However, $V$ is not maximal here.
We can see this by re-writing   $d v_1  = (u_1 + u_2)(u_1 + u_3)$,
and noting that only $2$ linearly independent elements of $U$
actually occur.  In other words, we can change the basis in $U$,
and move the `spare' generator of $U$ into $V$. Specifically,
define $\phi(u_1) = u_1$, $\phi(u_2) = u_2 + u_1$, $\phi(u_3) =
u_3 + u_1$ and $\phi(v_1) = v_1$, then extend to an algebra
automorphism of $\Lambda(U, V)$. If we set $d' = \phi^{-1}d\phi$,
then $\phi$ becomes an isomorphism of minimal models $\phi \colon
\Lambda(U, V; d') \to \Lambda(U, V; d)$. A simple check reveals
that $\Lambda(U, V; d')$ is two-stage with $d'(v_1) = u_2u_3$.
Thus we can write it as $\Lambda(U', V'; d')$, $U' = \langle u_2,
u_3 \rangle$ and $V' = \langle v_1, u_1 \rangle$.  In this latter
case, $V'$ is now of maximal dimension.
\end{example}

\begin{definition} Suppose $(\Lambda(U,V), d)$ is a two-stage minimal model.
We say that $V$ \emph{has maximal dimension}, or that the
two-stage decomposition displays $V$ with maximal dimension, if,
for any isomorphic two-stage minimal model $(\Lambda(U',V'), d' )
\cong (\Lambda(U,V), d)$, we have $\mathrm{dim}\,V' \leq
\mathrm{dim}\,V$.

Since we assume that the spaces of generators are
finite-dimensional, it is clear that every two-stage minimal model
has a decomposition that displays $V$ with maximal dimension.
Also, for the two-stage decomposition to display $V$ with maximal
dimension, then it is clearly necessary (but not sufficient) that
every generator from $U$ appear in some differential.
\end{definition}

We now give a consequence of the two-stage decomposition being
chosen so as to display $V$ of maximal dimension. This result does
not cover all two-stage cases, but is sufficient for our purposes.
We focus on the case in which the differential is quadratic.   For
parity of degree reasons, in this case we have $d(V) \subseteq
\Lambda^2(U^{\mathrm{odd}}) \oplus \Lambda^2(U^{\mathrm{even}})$.
Let $U^{\mathrm{odd}} = \langle u_1, \dots, u_p \rangle$,
$U^{\mathrm{even}} = \langle w_1, \dots, w_q \rangle$ and $V =
\langle v_1, \dots, v_r \rangle$ denote bases. From now on, we
will adopt this as standard notation, for the case of a two-stage
minimal model that may contain even-degree generators.

Since we assume the differential is quadratic,   for each basis
element $v_k$ of $V$, we may write
\begin{displaymath}
d v_k  = \sum_{1\leq i < j \leq p} a^k_{i, j}\,u_i u_j +
\sum_{1\leq i \leq j \leq q} b^k_{i, j}\,w_i w_j.
\end{displaymath}
Using the coefficients from the first of these sums, for
$k=1,\dots, r$, we define the skew-symmetric, $p \times p$ matrix
$M^k$,  by
\begin{displaymath}
M^k_{i,j}  =\left\{
\begin{array}{ll}
a^k_{i, j} & \text{if $i<j$}\\ - a^k_{i, j} & \text{if $i>j$}\\ 0
& \text{if $i=j$}\\
\end{array}\right.
\end{displaymath}
Now, let
$M=\bmatrix M^1\\
\vdots\\ M^r
\endbmatrix$ be the $pr \times p$ block marix formed by the $M^k$ as rows.

\begin{lemma}\label{lem:skew-symmetric span}
If $V$ has maximal dimension in  $\Lambda(U, V)$, then $\rank M =
p$. In particular, there is a $p \times rp$ matrix $N$ such that
$NM$ is the $p\times p$ identity matrix.
\end{lemma}

\begin{proof}
Let $U^*$ denote the dual  space of $U$, and for $u^*\in U^*$, let
$i_{u^*}$ denote the derivation of  $\Lambda(U, V)$ of degree $-1$
extending the linear map $i_{u^*}:U \oplus V\to \mathbb{Q}$
defined by $i_{u^*}(x+y)=u^*(x)$, for $x \in U$, $y \in V$. We
first show that, under our hypotheses, the Lie derivative
$\mathcal{L}\colon U^* \to \text{Hom}(V, U)$ defined by
$$
\mathcal{L}(u^*)=i_{u^*}d -(-1)^{|u^*|}di_{u^*}=i_{u^*}d
$$
is injective. (In fact, the injectivity of $\mathcal{L}$ is
equivalent to the maximality of $\dim V$.) Let  $K=\ker \mathcal
L$ and choose a  complement $X\subset U^*$ for $K$ so that $U^*=K
\oplus X$ and $U=K^{\perp} \oplus X^{\perp}$, where,  if $W\subset
U^*$, $W^{\perp}=\set{u\in U \mid \forall w \in  W,\,w(u)=0}$.
Then, as
$$\Lam U= (\Lam X^{\perp}
\otimes \Lam X^{\perp}) \oplus  (\Lam X^{\perp} \otimes \Lam^+
K^{\perp})\oplus (\Lam^+ K^{\perp} \otimes\Lam^+ K^{\perp}),$$ the
definition of $K$ shows that $dV \subset \Lam^+ K^{\perp}
\otimes\Lam^+ K^{\perp}$. Now define $U_1=K^{\perp}$ and
$V_1=V\oplus X^{\perp}$, with $d'_{U_1}=0$, $d_{|V}'=d_{|V}$ and
$d'_{X^{\perp}}=0$. Then $U_1\oplus V_1=U \oplus V$, and this
induces an isomorphism of   two-stage minimal models $\Lam
(U_1,V_1 ;d')\cong \Lam (U,V ;d)$. Since $V$ has maximal
dimension, $K$ must be 0.

Let $\set{u_1,\dots,u_p}$ be a basis of $U^{\odd}$, and
$\set{u_1^*,\dots,u_p^*}$ the dual basis. Then, as $V$ has maximal
dimension, we know that in particular the maps
$\mathcal{L}(u_j^*)=i_{u_j^*}d: V\to U^{\odd}$ for $j=1,\dots,p\,$
are linearly independent. In other words, if
$c=(c_1,\dots,c_p)^t\in \mathbb{Q}^p$, then $\sum_jc_j
i_{u_j^*}dv_k= \sum_{i,j} M^k_{i,j}u_ic_j =0$ for $k=1,\dots, r$,
implies that $c=0$. Because the $u_i$ are linearly independent,
this is equivalent to the statement $$\sum_{j} M^k_{i,j} c_j =0
\,\, \text{ for } i=1,\dots, p \text{ and }  k=1,\dots, r \quad
\Longrightarrow \quad c=0.$$ That is, $Mc=0 \implies c=0.$ Hence,
$\rank M= p.$  \end{proof}

To phrase our results, we use the following terminology:  We say
that a graded vector space $V$ is \emph{$n$-co-connected} if $V^i
= 0$ for $i \geq n$.

\begin{theorem}\label{thm:two stage even and odd}
Let $\Lambda(U,V)$ be an elliptic, two-stage minimal model. Assume
that $d V  \subseteq \Lambda^2U$, and that the two-stage
decomposition displays $V$ with maximal dimension.  Finally,
assume that $U^{\mathrm{odd}}$ and $U^{\mathrm{even}}$ satisfy one
of the following connectivity and co-connectivity hypotheses:
\begin{itemize}
\item[(A)]  $U^{\mathrm{odd}}$ is $(2r-1)$-connected and $U^{\mathrm{even}}$
is $(2r+2)$-co-connected.
\item[(B)]   $U^{\mathrm{odd}}$ is $(2s+1)$-co-connected and
  $U^{\mathrm{even}}$ is $(4s-4)$-connected.
\end{itemize}
Then $\mathrm{rk}_0\big(\Lambda(U, V)\big) \leq \mathrm{dim}\,V -
\mathrm{dim}\,U^{\mathrm{even}}$.
\end{theorem}

\begin{proof}
If $\mathrm{rk}_0\big(\Lambda(U, V)\big) = n$, then we have a
K-S-extension $(\Lambda A_n \otimes \Lambda(U, V), D)$ as in
(\ref{eq:Borel fibration}) that has finite-dimensional cohomology.
We claim that in any such minimal model,  we can assume that
$D(U^{\mathrm{odd}})$ is contained in $I(U^{\mathrm{odd}}\oplus
V)$, the ideal in $\Lambda A_n \otimes \Lambda(U, V)$ generated by
$U^{\mathrm{odd}}\oplus V$.

Allowing this claim for the time-being, we appeal to some results
of Halperin concerning elliptic minimal models. First of all, to
any elliptic minimal model $(\Lambda W, d)$, there is an
associated \emph{pure model}, denoted $(\Lambda W, d_{\sigma})$,
which is defined by adjusting the differential $d$ to $d_{\sigma}$
as follows:  We set $d_{\sigma} = 0$ on each even degree generator
of $W$, and on each odd degree generator $w \in W$, we set
$d_{\sigma}(w)$ equal to the part of $d(w)$ contained in
$\Lambda(W^{\mathrm{even}})$.  One checks that this defines a
differential $d_{\sigma}$ on $\Lambda W$, and thus we obtain a
pure (minimal) model $(\Lambda W, d_{\sigma})$. Then, by
\cite[Prop.1]{Halperin-finiteness} (cf.~also
\cite[Prop.32.4]{F-H-T00}), $\mathrm{dim}\,H(\Lambda W, d)$ is
finite-dimensional if and only if $\mathrm{dim}\,H(\Lambda W,
d_{\sigma})$ is finite-dimensional. Applying all this to the
minimal model $(\Lambda A_n \otimes \Lambda(U, V), D)$, we obtain
the following:

 From the claim, $(\Lambda A_n \otimes \Lambda(U, V),
D_{\sigma})$ satisfies $D_{\sigma}(U^{\mathrm{odd}}) = 0$, and
therefore $(\Lambda A_n \otimes \Lambda(U, V), D_{\sigma}) \cong
(\Lambda U^{\mathrm{odd}}, D_{\sigma} = 0) \otimes (\Lambda A_n
\otimes \Lambda(U^{\mathrm{even}}, V), D_{\sigma})$.  Since
$(\Lambda A_n \otimes \Lambda(U, V), D_{\sigma})$ has
finite-dimensional cohomology, so does $(\Lambda A_n  \otimes
\Lambda(U^{\mathrm{even}}, V), D_{\sigma})$.  It follows that
$(\Lambda A_n  \otimes \Lambda(U^{\mathrm{even}}, V), D_{\sigma})$
is an elliptic minimal model, and therefore that its homotopy
Euler characteristic is non-positive
\cite[Th.1]{Halperin-finiteness} (cf.~also
\cite[Prop.32.10]{F-H-T00}). This implies that $n \leq
\mathrm{dim}\,V - \mathrm{dim}\, U^{\mathrm{even}}$, as required.

It only remains to establish the claim.  We do this by a careful
analysis of the terms that can occur in the differentials. Our
argument is essentially the same as that in
\cite[Th.4.6]{Allday-Puppe85b}. First we write, for each $v_k$ in
the basis of $V$,
\begin{equation}
\label{eq:differential D} D v_k  = \sum_{1\leq i < j \leq p}
a^k_{i, j}\,u_i u_j \ +\  \sum_{1\leq i \leq j \leq q} b^k_{i,
j}\,w_i w_j\  +\   \text{terms in the ideal $I(A_n)$}.
\end{equation}
The coefficients $a^k_{i, j}$ and $b^k_{i, j}$ are scalars.   Our
strategy in either case (A) or case (B) is the same.  We apply $D$
to (\ref{eq:differential D}), obtaining $D^2 v_k  = 0$ on the
left-hand side.  By focussing on certain terms present on the
right-hand side, we obtain sufficient information to complete our
argument.

Consider hypothesis (A) first. Here, from the (co-)connectivity
restrictions, we see that $D(U^{\mathrm{even}}) \subseteq
\Lambda^{+}A_n\otimes\Lambda  U^{\mathrm{even}} \otimes
\Lambda^{+}V$.  Therefore, display the terms in
(\ref{eq:differential D}) as follows. Write
\begin{equation}
\label{eq:differential D(A)} \begin{aligned} D v_k  = \ &
\sum_{1\leq i < j \leq p} a^k_{i, j}\,u_i u_j +\
   \sum_{1\leq i < j \leq p} P^k_{i,j}\,u_i u_j \  +\
\sum_{\stackrel{\scriptstyle{i=1,\dots,r}}{j=1,\dots,p}} R^k_{i,
j}\,v_iu_j
\\ &\  + \ P_k \   +\  \text{terms in the ideal
$I\big(\Lambda^{2}(V)\oplus\Lambda^{\geq
4}(U^{\mathrm{odd}},V)\big)$}.
\end{aligned}
\end{equation}
Here, $P_k$ includes all terms from $\Lambda A_n \otimes
\Lambda(U^{\mathrm{even}})$ present in (\ref{eq:differential D}),
and the coefficients $P^k_{i, j}$ and $R^k_{i, j}$ denote terms in
$\Lambda^{+} A_n$. After applying $D$ to this, we will consider
only contributions in $\big(\Lambda^{+} A_n\otimes
\Lambda(U^{\mathrm{even}})\big)\cdot U^{\mathrm{odd}}$. Define a
derivation $\delta$ of $\Lambda A_n \otimes \Lambda (U, V)$ as
follows. On generators from $U^{\mathrm{odd}}\oplus V$, set
$\delta$ equal to the part of $D \colon U^{\mathrm{odd}}\oplus V
\to \Lambda A_n\otimes \Lambda(U, V)$ whose image is contained in
$\Lambda^{+} A_n\otimes \Lambda(U^{\mathrm{even}})$. On the
remaining generators, set $\delta(U^{\mathrm{even}}) = 0$. Extend
$\delta$ to a derivation on $\Lambda(A_n, U, V)$. Applying $D$ to
(\ref{eq:differential D(A)}), and collecting terms in
$\big(\Lambda^{+} A_n\otimes \Lambda(U^{\mathrm{even}})\big)\cdot
U^{\mathrm{odd}}$, we have
\begin{displaymath}
0 = \delta(\sum_{1\leq i < j \leq p} a^k_{i, j}\,u_i u_j) +
\delta(\sum_{1\leq i < j \leq p} P^k_{i, j}\,u_i u_j) +
\sum_{\stackrel{\scriptstyle{i=1,\dots,r}}{j=1,\dots,p}}
\delta(R^k_{i, j}v_i)u_j .
\end{displaymath}
We fix an index $l$, with $1 \leq l \leq p$, and collect all
coefficients of $u_l$ from this equation, to obtain
\begin{displaymath}
0 = \sum_{i < l} a^k_{i, l}\,\delta(u_i) - \sum_{l < j} a^k_{l,
j}\,\delta(u_j) + \sum_{i<l} \delta(P^k_{i, l}\,u_i) - \sum_{l<j}
\delta(P^k_{l, j}\,u_j) + \sum_{i=1,\dots,r} \delta(R^k_{i,l}v_i).
\end{displaymath}
Since $P^k_{i, j}, R^k_{i, j} \in \Lambda^{+} A_n$, and since
$M^k$ is skew-symmetric, we obtain
\begin{displaymath}
\sum_{i \not= l} a^k_{l,i}\,\delta(u_i)  =  \delta(\chi_{k, l})
\end{displaymath}
for some element $\chi_{k, l}\in \Lambda^{+}A_n\cdot(U\oplus V)$.
Let $\mathbf{u}=(u_1,\dots,  u_p)^t$, $\sigma_k=(\chi_{k,
1},\dots,\chi_{k, p})$, and $\mathbf{\chi}=(\sigma_1
,\dots,\sigma_r)^t$, and let $\delta$ act component-wise. Then the
previous equation can be restated as
$$\delta(M\mathbf{u})=\delta \mathbf{\chi}.
$$

  Now we
apply \lemref{lem:skew-symmetric span} to obtain a $p \times rp$
matrix $N$ such that $$\delta\mathbf{u}
=\delta(NM\mathbf{u})=\delta (N\mathbf{\chi}).
$$   Finally, we observe that this implies for
each basis element of $U^{\mathrm{odd}}$, we have $\delta(u_i) =
\delta(\chi_i)$ for some decomposable $\chi_i \in
\Lambda^{+}A_n\cdot(U^{\mathrm{odd}}\oplus V)$. Returning to the
full differential $D$, we rephrase this as follows: For each basis
element $u_i$ of $U^{\mathrm{odd}}$, there is some $\chi_i \in
\Lambda^{+}A_n\cdot(U^{\mathrm{odd}}\oplus V)$ for which $D(u_i -
\chi_i) \in I(U^{\mathrm{odd}}, V)$.

As in \exref{ex:non-maximal V}, we now change the basis in
$U^{\mathrm{odd}}$. Define an isomorphism of algebras $\phi \colon
\Lambda A_n \otimes \Lambda(U, V) \to \Lambda A_n \otimes
\Lambda(U, V)$ on generators by setting $\phi(u_i) = u_i - \chi_i$
for each basis element $u_i$ of $U^{\mathrm{odd}}$, and $\phi =
\text{id}$ on generators of $A_n$, $U^{\mathrm{even}}$ and $V$.
Then extend $\phi$ to an algebra isomorphism.  Now define a
differential by $D' = \phi^{-1}D\phi$, so that $\phi$ becomes an
isomorphism of minimal models $\phi \colon (\Lambda A_n \otimes
\Lambda(U, V), D') \to (\Lambda A_n \otimes \Lambda(U, V), D)$. An
easy check now confirms that $D'(A_n) = 0$ and
$D'(U^{\mathrm{odd}}) \subseteq I\big(U^{\mathrm{odd}}, V \big)$,
as desired.

A similar argument is used under hypothesis (B). Here, if some
coefficient $a^k_{i, j}\not= 0$ in (\ref{eq:differential D}), then
the degree of $D v^k $ must be at most $4s-2$.  But under
hypothesis (B), terms in $\Lambda^2(U^{\mathrm{even}})$ start in
degree $8s-4$ and terms in $\Lambda^{+}A_n\cdot\Lambda
U^{\mathrm{even}}$ start in degree $4s$. In this case, therefore,
whenever some $a^k_{i, j}\not=0$, then the remaining terms of
(\ref{eq:differential D}) do not include any in the ideal
$I(U^{\mathrm{even}})$. So suppose that $v^k$ is of suitable
degree for $D(v^k)$ to contain non-zero terms $a^k_{i, j}\,u_i
u_j$. Display the terms in (\ref{eq:differential D}) as follows.
Write
\begin{equation}
\label{eq:differential D(B)} \begin{aligned} D(v_k) = \sum_{1\leq
i < j \leq p} a^k_{i, j}\,u_i u_j \ &+\
   \sum_{1\leq i < j \leq p} P^k_{i,
j}\,u_i u_j  +
\sum_{\stackrel{\scriptstyle{i=1,\dots,r}}{j=1,\dots,p}} R^k_{i,
j}\,v_iu_j   \\ & + P_k + \text{terms in $I\big( \Lambda^{\geq2}V
\oplus \Lambda^{\geq 4}U^{\mathrm{odd}}\big)$},
\end{aligned}
\end{equation}
where $P_k$ denotes a term from $\Lambda^{+}A_n$  and the
coefficients $P^k_{i, j}$ and $R^k_{i, j}$ denote terms in
$\Lambda^{+} A_n$. After applying $D$ to this, we will consider
only contributions in $(\Lambda^{+} A_n)\cdot U^{\mathrm{odd}}$.
To this end, define a derivation $\delta$ of $\Lambda A_n \otimes
\Lambda (U, V)$ as follows. On generators from
$U^{\mathrm{odd}}\oplus V$, set $\delta$ equal to the part of $D
\colon U^{\mathrm{odd}}\oplus V \to \Lambda A_n\otimes \Lambda(U,
V)$ whose image is contained in $\Lambda^{+}A_n$. On the remaining
generators, set $\delta(U^{\mathrm{even}}) = 0$. Applying $D$ to
(\ref{eq:differential D(B)}), and collecting terms in
$(\Lambda^{+} A_n)\cdot U^{\mathrm{odd}}$ yields
\begin{displaymath}
0 = \delta(\sum_{1\leq i < j \leq p} a^k_{i, j}\,u_i u_j) +
\delta(\sum_{1\leq i < j \leq p} P^k_{i, j}\,u_i u_j) +
\sum_{\stackrel{\scriptstyle{i=1,\dots,r}}{j=1,\dots,p}}
\delta(R^k_{i, j}v_i)u_j ,
\end{displaymath}
for each index $k$ that has at least one coefficient $a^k_{i, j}$
non-zero. From here, the argument proceeds exactly as that for
hypothesis (A).

This completes the proof of the claim, under either hypothesis. By
the first part of the proof, the result follows.
\end{proof}

There are other situations in which a similar approach gives an
upper bound on the toral rank. We give one more result that moves
away from the quadratic differential restriction.

\begin{theorem}\label{thm:stable and separated}
Let $\Lambda(U, V)$ be a two-stage minimal model with odd degree
generators only and assume that the two-stage decomposition
displays $V$ with maximal dimension. Suppose there are positive
integers $r$ and $s$, with $r \leq s \leq 2r$, for which $U =
U^{\mathrm{odd}}$ is $(r-1)$-connected and $(s+1)$-co-connected.
Suppose further that there are positive integers $t$ and $u$, with
$s \leq t \leq u \leq s + r$, for which $V$ is $(t-1)$-connected
and $(u+1)$-co-connected.  Then $\mathrm{rk}_0\big(\Lambda(U, V;
d)\big) \leq \mathrm{dim}\,V$,
\end{theorem}

\begin{proof}
If $\mathrm{rk}_0\big(\Lambda(U, V; d)\big) = n$, then we have an
extension $(\Lambda A_n \otimes\Lambda(U, V), D)$ with
finite-dimensional cohomology. In such an extension, from the
assumptions on the degrees of the generators, we must have $D u_i
\in \Lambda A_n$ for each generator $u_i$ of $U$. We claim that in
fact $D u_i  = 0$. Suppose the shortest length of any non-zero
term that occurs in some polynomial $D u_i $ is $l$. Then as in
the proof of \thmref{thm:two stage even and odd}, define a
derivation $\delta$ on $U$ by setting $\delta u_i $ equal to the
term of length $l$ that occurs in $D u_i $, or zero if there is no
such term. Define $\delta$ to be zero on all generators of $A_n$
and $V$, and extend as a derivation to $\Lambda A_n
\otimes\Lambda(U, V)$. Without loss of generality, we can suppose
that the non-zero terms $\delta u_i $ are linearly independent.
Otherwise, we can change basis within $U$ to make these so.  Next,
consider $D v $ for a generator in $V$.  We have $D v  = dv +
\chi$, for $dv \in \Lambda U$ and some $\chi$ in the ideal of
$\Lambda A_n \otimes\Lambda(U, V)$ generated by $A_n$. The
(co)-connectivity hypotheses on $V$ imply that $D v  \in \Lambda
A_n \otimes\Lambda U$, for any $v \in V$.  Therefore, $\chi \in
\Lambda^{+} A_n \otimes\Lambda U$ and hence $D(\chi) \in
\Lambda^{>l} A_n \otimes\Lambda U$.   Consequently, when we equate
terms in $\Lambda^{l} A_n \otimes\Lambda U$ that arise in the
equation $0 = D^2 v  = D(dv) + D(\chi)$, we have $0 = \delta(dv)$.
  By the assumption that non-zero $\delta u_i $'s are
linearly independent, it follows from an easy argument that
$\delta u_i = 0$ for each $u_i$ that occurs in the differential
$dv$.  Finally, our assumption that $V$ is taken as large as
possible implies that each generator $u_i$ does occur in at least
one differential $dv$.  This implies that the shortest length
terms $\delta u_i $ are zero, and an induction completes the proof
of our claim, namely that $D u_i  = 0$ for each $u_i \in U$.

So far, we have argued that in any extension $(\Lambda A_n
\otimes\Lambda(U, V), D)$, our hypotheses imply that $D(U) = 0$.
Now an argument as in the second paragraph of the proof of
\thmref{thm:two stage even and odd} shows that $n \leq
\mathrm{dim}\,V$.
\end{proof}

\begin{corollary}\label{cor:two stage even and odd}
Let $X$ be a simply connected, elliptic space with two-stage
minimal model.  If the minimal model satisfies either the
hypotheses of \thmref{thm:two stage even and odd}, or those of
\thmref{thm:stable and separated}, then \conjref{conj:TRC} holds
for $X$.
\end{corollary}

\begin{proof}
Combine \thmref{thm:general 2-stage} with the results mentioned in
the statement.
\end{proof}

It is natural to ask whether $\mathrm{rk}_0(X) \leq
\mathrm{dim}\,V - \mathrm{dim}\,U^{\mathrm{even}}$ for a general
two-stage, elliptic space.  This is not the case in general, as we
illustrate in \exref{ex:rk>dimV - dimUev} below. Indeed, under
hypotheses that include all cases in which the two-stage minimal
model has odd-degree generators only, we can actually reverse the
direction in this inequality.

\begin{lemma}\label{lem:upper bound}
  Suppose $\Lambda  U \otimes \Lambda  V$
is an elliptic two-stage minimal model with
$U^{\mathrm{even}}=U^{2n}$ for some fixed $n\in \mathbb{N}$. Then,
$\mathrm{rk}_0(\Lambda  U \otimes \Lambda V)\geq \dim V-\dim
U^{\mathrm{even}}$.
\end{lemma}

\begin{proof} Consider the extension sequence
\begin{displaymath}
(\Lambda  U^{\mathrm{odd}},0) \to  (\Lambda
U^{\mathrm{odd}}\otimes \Lambda  U^{\mathrm{even}} \otimes \Lambda
V,d) \to (\Lambda U^{\mathrm{even}} \otimes \Lambda  V, \bar d).
\end{displaymath}
Since $V=V^{\mathrm{odd}}$, the right-hand term is now  pure (and
elliptic) with $U^{\mathrm{even}}=U^{2n}$.  By
\cite[Lem.3.3]{Jess}, there is an isomorphism of two-stage models
$$\rho \colon  (\Lambda U^{\mathrm{even}} \otimes \Lambda V, \bar d) \to
(\Lambda U^{\mathrm{even}} \otimes \Lambda  (V' \oplus V''), \bar
d')$$ such that

\begin{enumerate}
\item  $\dim H(\Lambda  U^{\mathrm{even}} \otimes \Lambda  V',\bar d) <
\infty$,
\item  $\dim V'=\dim U^{\mathrm{even}}$,
\item $\rho$ is the identity on $U^{\mathrm{even}}$, and
\item $\rho(V) \subset V'\oplus
V'' \oplus (\Lambda ^{+} U^{\mathrm{even}}\otimes (V'\oplus V''))
$.
\end{enumerate}

These conditions imply that $\rho$ induces an isomorphism $V \cong
V'\oplus V'',$ and that $\rho^{-1}(V'\oplus V'') \subset V\oplus
(\Lambda ^{+} U^{\mathrm{even}}\otimes V ) $. We now extend $\rho$
to an isomorphism of algebras $$\tilde \rho \colon \Lambda  U
\otimes\Lambda  V \to \Lambda  U \otimes \Lambda (V' \oplus V'')$$
    by letting it   be the identity on $\Lambda  U$. Define a
derivation $d'$ on $\Lambda  U \otimes \Lambda  (V' \oplus V'')$
by setting $d'U=0$, and $ d' v'= d\rho^{-1}(v'), $ for $v'\in (V'
\oplus V'')$.  In order to show that $$\tilde \rho \colon (\Lambda
U \otimes\Lambda V ,d) \to (\Lambda  U \otimes \Lambda  (V' \oplus
V''),d')$$ is an isomorphism of models, it suffices to show that
$(d')^2v'=0$ for $v'\in V' \oplus V'',$ so we compute:
\begin{align*}
(d')^2v' &=d'(d\rho^{-1}(v'))\\
   &\in d'(d( V\oplus (\Lambda ^{+} U^{\mathrm{even}}\otimes V )) \\
   &\subset d'(\Lambda  U) \\
   &=\{ 0\}.
\end{align*}
Thus, we may assume that $V=V^{\mathrm{odd}}=V'\oplus V''$, with
$\dim H(\Lambda  U \otimes \Lambda  V') < \infty$, and $\dim
V'=\dim U^{\mathrm{even}}$.  Now suppose $V'' = \langle v_1,
\dots, v_n \rangle$, and define a K-S-extension
$$(\Lambda A_n,0) \to (\Lambda A_n \otimes \Lambda (U, V), D) \to (\Lambda
(U, V), d),
$$
in which $A_n = \langle a_1, \dots, a_n \rangle$ with$|a_i|=2$, by
$D v_i  = d v_i  + a_i^{\frac{|v_i| + 1}{2}}$, and $d=D$
otherwise.  A standard argument (using  the associated pure model,
as in the proof of Theorem \ref{thm:two stage even and odd}), now
shows that $\mathrm{dim}\,H(\Lambda A_n \otimes \Lambda (U, V),
D)$ is finite-dimensional, so that $\mathrm{rk}_0\Lambda (U, V)
\ge n=\dim V''=\dim V-\dim U^{\mathrm{even}}$.
\end{proof}

Consequently, we obtain values for the rational toral rank in the
following cases.

\begin{corollary}\label{cor:exact rank}
Let $X$ be a simply connected elliptic space with two-stage
minimal model $\Lambda (U, V)$.  Suppose that the two-stage
decomposition displays $V$ with maximal dimension.
\begin{itemize}
\item[(\ref{cor:exact rank}.1)] If the minimal model has quadratic
differential, satisfies one of the (co-) connectivity conditions
of \thmref{thm:two stage even and odd} and also satisfies
$U^{\mathrm{even}}=U^{2n}$ for some fixed $n$, then
$\mathrm{rk}_0(X) = \dim V-\dim U^{\mathrm{even}}$.
\item[(\ref{cor:exact rank}.2)] If the minimal model has odd-degree
generators only, and satisfies the (co-) connectivity conditions
of \thmref{thm:stable and separated}, then $\mathrm{rk}_0(X) =
\dim V$.
\end{itemize}
\end{corollary}

\begin{proof}
Combine \lemref{lem:upper bound} with the results mentioned in the
statement.
\end{proof}

\begin{remarks}\label{rems:special cases}
We can specialize (\ref{cor:exact rank}.1) in a number of
interesting directions.  For example, if $U^{\mathrm{even}}= 0$,
then we retrieve the identification $\mathrm{rk}_0(X) = \dim V$ of
\thmref{thm:two stage quadratic}. As a second example, we can
restrict to the case in which $U^{\mathrm{odd}}= 0$.  This is the
case in which the minimal model is \emph{pure} with quadratic
differential, and also satisfies $U^{\mathrm{even}}=U^{2n}$ for
some fixed $n$.  Here, we obtain $\mathrm{rk}_0(X) = -
\chi_{\pi}(X)$.
\end{remarks}

\section{Examples, Comments and Questions}

In this section we mention various examples and results. Our focus
here is more on the exact toral rank, rather than the bound of
\conjref{conj:TRC}.

\subsection{The relation between $\mathrm{rk}_0(X)$ and $\mathrm{dim}\,V -
\mathrm{dim}\,U^{\mathrm{even}}$} We begin with the simplest
example that we can find to illustrate that the inequality
$\mathrm{rk}_0(X) \leq \mathrm{dim}\,V -
\mathrm{dim}\,U^{\mathrm{even}}$ does not hold in general for a
two-stage, elliptic space with both even and odd generators.

\begin{example}\label{ex:rk>dimV - dimUev}
Consider the two-stage minimal model $\Lambda(u_1, u_2, u_3, u_4,
u_5, u_6, w, v)$, with $|u_1| = |u_2| = |u_3| = 3$, $|u_4| = 7$,
$|u_5| = 23$, $|u_6| = 27$, $|w| = 18$ and $|v| = 35$, and the
single non-trivial differential given by $d v  = w^2 - u_1 u_2 u_4
u_5 - u_1 u_2 u_3 u_6$.   The associated pure model satisfies
$(\Lambda(U, V), d_{\sigma}) \cong (\Lambda(u_1, u_2, u_3, u_4,
u_5, u_6), d_{\sigma} = 0) \otimes (\Lambda(w, v), d_{\sigma})$,
with $d_{\sigma}(w) = 0$ and $d_{\sigma}(v) = w^2$.  Now
$H(\Lambda(w, v), d_{\sigma}) \cong \Lambda(w)/(w^2)$, so  $
(\Lambda(U, V), d)$ is elliptic.

The two-stage decomposition $U = \langle u_1, u_2, u_3, u_4,$
$u_5, u_6, w \rangle$ and $V = \langle v \rangle$ displays $V$
with maximal dimension, so we have $\mathrm{dim}\,V -
\mathrm{dim}\,U^{\mathrm{even}} = 0$.  We now show that
$\mathrm{rk}_0\big(\Lambda(U, V)\big)) \geq  1$. Let $a$ be a
generator of degree $2$, and define a differential $D$ on $\Lambda
a)\otimes\Lambda(U, V)$ as follows: $D a  = D  u_1  = D u_2  = D
u_3  = D u_4  = D w  = 0$, $D u_5  = a^3w$, $D u_6  = a^2wu_1u_2$
and $D v  = w^2 - u_1 u_2 u_4 u_5 - u_1 u_2 u_3 u_6 + a^{18} -
au_6u_4$.  A straightforward check shows that this defines a
differential.  We show that $\big(\Lambda a \otimes\Lambda(U, V),
D\big)$ has finite-dimensional cohomology.  The associated pure
model in this case is $(\Lambda a \otimes\Lambda(U, V),
D_{\sigma}) \cong (\Lambda(u_1, u_2, u_3, u_4, u_6), D_{\sigma} =
0) \otimes (\Lambda(a, w, u_5, v), D_{\sigma})$, with $D_{\sigma}
a  = 0$, $D_{\sigma}(w) = 0$, $D_{\sigma} u_5  = a^3w$ and
$D_{\sigma} v  = w^2+a^{18}$.   We observe that $D_{\sigma}(a^3v)
= a^3w^2 +a^{21} = D_{\sigma}(u_5w) + a^{21}$, hence
$D_{\sigma}(a^3v - u_5w) = a^{21}$, so
\cite[Prop.1]{Halperin-finiteness} shows that $H(\Lambda a
\otimes\Lambda(U, V), D_{\sigma})$ is finite-dimensional.  It
follows that $H(\Lambda a \otimes\Lambda(U, V), D)$ is finite
dimensional. Thus $\mathrm{rk}_0\big(\Lambda(U, V)\big)) \geq  1$.
\end{example}

We assert that more work will show $\mathrm{rk}_0\big(\Lambda(U,
V)\big)) = 1$ in \exref{ex:rk>dimV - dimUev}.  We also remark that
whilst \conjref{conj:TRC} does not follow from \thmref{thm:general
2-stage} in this example, it is nonetheless easily confirmed here.

Next, we specialize to the case of odd generators only.
Unfortunately, this restriction alone does not give us the
inequality $\mathrm{rk}_0(X) \leq \mathrm{dim}\,V$.  The following
example  lies immediately outside the hypotheses of
\thmref{thm:stable and separated}, at least as far the generators
of $U$ are concerned. We thank Yves F{\'e}lix for  this example.

\begin{example}\label{ex:two-stage odd}
Let $\Lambda W = \Lambda(U, V)= \Lambda(u_1, u_2, u_3, u_4, u_5,
v_1, v_2)$ be a two-stage minimal model with degrees and
differential as follows: $|u_1| = |u_2| = |u_3| = |u_4| = 3$,
$|u_5| = 7$, $|v_1| = 9$, $|v_2| = 11$, $d u_i  = 0$ for each $i$,
$d v_1  = u_1 u_5$ and $d v_2  = u_1 u_2 u_3 u_4$. Then
$\mathrm{dim}\,V = 2$. We show that $\mathrm{rk}_0\big(\Lambda
W\big) \geq 3$.  First we describe an extension $\Lambda A
\otimes\Lambda W, D$, with $A = \langle a_1, a_2, a_3\rangle$. Set
$D = 0$ on $\{ a_1, a_2, a_3, u_1, u_2, u_3 \}$, and $D u_4  =
-a_1^2$, $D u_5  = a_1 u_1 u_2$, $D v_1  = u_1 u_5 + a_2^5$ and $D
v_2  = u_1 u_2 u_3 u_4 + u_3 u_5 a_1 + a_3^6$. A straightforward
check shows that $D$ is a differential. We argue exactly as in the
latter part of \exref{ex:rk>dimV - dimUev} to show that $H(\Lambda
A \otimes\Lambda W, D)$ is finite-dimensional. It follows  that
$\mathrm{rk}_0\big(\Lambda W\big) \geq 3$.
\end{example}

Once again, \conjref{conj:TRC} does hold for \exref{ex:two-stage
odd}, although we are unable to conclude this from our preceding
results.  To see this, use \cite[Th.4.6]{Allday-Puppe85b} to
conclude $\mathrm{rk}_0\big(\Lambda W\big) \leq 5$.  (In
\exref{ex:two-stage odd}, the differential $d v_1  = u_1 u_5$
means that $u_1$ and $u_5$ correspond to non-central generators.)
Moroever, the fibration with $u_1$ as base has fibre $ \Lambda(
u_2, u_3, u_4, u_5, v_1, v_2) $ with zero differential, and
$\Lambda( u_2, u_3, u_4, u_5)$ is clearly in the kernel of the
associated Wang derivation $\theta^*$. As in the proof of Theorem
\ref{prop:odd-free}, we conclude that $\dim H(\Lam W)=2 \dim \ker
\theta^*\ge 2. 2^4=2^5  \geq 2^{\mathrm{rk}_0(\Lambda W)}$.

Another example in which $\mathrm{rk}_0(X) > \mathrm{dim}\,V -
\mathrm{dim}\,U^{\mathrm{even}}$ is given in  \exref{ex:Steve's
example}. On the other hand, we currently have no example of a two
stage minimal model for which $\mathrm{rk}_0(X) < \mathrm{dim}\,V
- \mathrm{dim}\,U^{\mathrm{even}}$. This leads us to wonder
whether the hypotheses of \lemref{lem:upper bound} might be
relaxed considerably. Note that in certain cases, the inequality
$\mathrm{rk}_0(X) \geq \mathrm{dim}\,V -
\mathrm{dim}\,U^{\mathrm{even}}$ can be combined with other
information to identify the toral rank,
  and we have seen instances of this in \corref{cor:exact rank}. As another
example, whenever this latter inequality holds in a pure case, it
identifies the toral rank as equal to $-\chi_{\pi}$
(cf.~\remsref{rems:special cases}).

\subsection{Products} Now we consider the question of how the toral rank
behaves with respect to products. Although our previous results
concerned the two-stage case, the comments here are not restricted
to that case. It is easy to see that, in general, we have the
inequality $\mathrm{rk}_0(X\times Y) \geq \mathrm{rk}_0(X) +
\mathrm{rk}_0(Y)$. The following example shows that we may
sometimes have strict inequality.

\begin{example}
\label{ex:Steve's example} Consider two-stage, elliptic minimal
models $\mathcal{M}  = \Lambda(x, y)$ with $|x| =  12$, $|y| =
23$, and $d y  = x^2$; and $\mathcal{N}  = \Lambda(u_1, u_2, u_3,
u_4, u_5, u_6, w, v)$, with $|u_1| = |u_2| = |u_3| = |u_4| =  3$,
$|u_5| =  5$, $|u_6| =  19$, $|w| =  18$, $|v| =  35$, and $d v  =
w^2 + u_1 u_2 u_3 u_4 u_5 u_6$, all other differentials being
zero.

   The homotopy
Euler characteristic  bound yields $\mathrm{rk}_0(\mathcal{M} ) =
0$. We will show that $\mathrm{rk}_0(\mathcal{N} ) = 0$ by arguing
that any extension of the form $(\Lambda a \otimes\mathcal{N} ,
D)$ cannot have finite-dimensional cohomology.

Let $U' = \langle u_1,\dots,u_5 \rangle$  and denote $
\Lambda^{\set{i_1,\dots, i_k}}U' :=\bigoplus_{m=1}^k\Lambda^{i_m}
U'$.  For degree reasons, any such $D$ satisfies $D u_i  =
\alpha_i\, a^{ 2}$, for $i = 1,\dots, 4$, $D u_5  = \alpha_5\, a^{
3}$, $D u_6  = \alpha_6\, a^{ 10} + \mu\, a w + \Phi$, and
\begin{equation}
\label{eq:differential Steve's} D v  = w^2 + u_1 u_2 u_3 u_4 u_5
u_6 + w \Psi + u_6 \Gamma + \Omega,
\end{equation}
where the $ \alpha_i$ and $\mu$ are scalars, $\Phi,Dw  \in
\Lambda^+a\otimes \Lambda^+U', \,\Gamma \in \Lambda^{\ge
3}a\otimes \Lambda^{\set{1,3}}U' $ and $ \Psi, \Omega
\in\Lambda^{+}a\otimes \Lambda U'$.

  Our
strategy is to show that each $\alpha_i = 0$, and then to follow a
similar argument as in previous examples.   Applying $D$ to
(\ref{eq:differential Steve's}) and equating terms that contain
$u_6$, we find
$$0=\sum_{i=1}^4 \pm \alpha_ia^{ 2}u_1\dots \hat{u_i}\dots u_4 \,u_5
+u_1\dots u_4  \alpha_5a^{ 3}  - D\Gamma.$$
Since $D\Gamma \in \Lambda^{\ge  5}a\otimes \Lambda U'$,   we
immediately find $\alpha_i=0$ for $i=0,\dots ,5$, and so   $D
U'=0$. Now,  $ 0 =D^2u_6 =\mu a Dw $. So either $\mu = 0$ or $Dw =
0$.   Moreover, the $w$ component of  $ D^2v$ being zero implies
that
$$0=2 Dw  -\mu  a u_1\dots u_5 + \mu a \Gamma.$$
If $\mu = 0$, then this equation implies $Dw = 0$.  On the other
hand, suppose $\mu \not= 0$.  Since $\Gamma \in \Lambda^{\ge
3}a\otimes \Lambda^{\set{1,3}}U $, by multiplying this last
equation by $\mu a$, we conclude that in fact $\mu=0$, and hence
that $Dw=0$.

The equation $D^2v=0 $  now implies
\begin{equation}
\label{D2v}0=  -\alpha_6 a^{ 10}u_1\dots u_5    +(\alpha_6a^{ 10}
+ \Phi)
   \Gamma.
\end{equation}
Suppose $\alpha_6\not=0$. Then, upon considering the component of
(\ref{D2v}) in $\Lam ^+a\otimes \Lam ^1 U'$, we find that
$\Gamma\in \Lam ^+a\otimes\Lam^{3}U$. Since $\Phi \in \Lam
^+a\otimes\Lam^{\ge 2} U'$,  this shows that $\Phi \Gamma \in \Lam
^+a\otimes\Lam^{5} U'$. Thus, (\ref{D2v}) now yields $\Gamma=0$
and so $\a_6=0$ in any case.

The  argument so far has shown $\mu=\alpha_i = 0$ for $i = 1,
\ldots, 6$.  Therefore, for any minimal model $(\Lambda a
\otimes\mathcal{N} , D)$, the corresponding pure model satisfies
$(\Lambda a \otimes\mathcal{N} , D_{\sigma}) \cong (\Lambda(u_1,
u_2, u_3, u_4, u_5, u_6),0)\otimes(\Lambda(a, w, v), D_{\sigma})$.
Since $(\Lambda(a, w, v), D_{\sigma})$ has homotopy Euler
characteristic $\chi_{\pi} = +1$, it cannot have
finite-dimensional cohomology, and hence neither can $(\Lambda a
\otimes\mathcal{N} , D_{\sigma})$.  As before, it follows that
$(\Lambda a \otimes\mathcal{N} , D)$ does not have
finite-dimensional cohomology, and hence that
$\mathrm{rk}_0(\mathcal{N} ) =0$.

Now consider the product $\mathcal{M} \otimes \mathcal{N} $. We
show that $\mathrm{rk}_0(\mathcal{M} \otimes \mathcal{N} ) \geq
1$, by displaying an extension $(\Lambda a \otimes \mathcal{M}
\otimes \mathcal{N} , D)$, in which $H(\Lambda a \otimes
\mathcal{M}\otimes \mathcal{N}, D)$ is finite-dimensional. The
differential is as follows: $Dx = u_1 u_2 u_3 a^2$, $Dy = x^2 + 2
u_1 u_6 a$, $Du_1 = Du_2 = Du_3 = Du_4 = 0$, $Du_5 = a^{ 3}$,
$Du_6 = u_2 u_3 x a$, $Dw = 0$ and $Dv = w^2 + u_1 u_2 u_3 u_4 u_5
u_6 + xa u_4u_6 + ya^{ 2}u_2u_3u_4$. A careful check reveals that
this defines a differential.  The associated pure model is now
$(\Lambda a \otimes \mathcal{M}  \otimes \mathcal{N} ,
D_{\sigma})$ with the only non-trivial differentials
$D_{\sigma}(y) = x^2$, $D_{\sigma}(u_5) = a^{ 3}$ and
$D_{\sigma}(v) = w^2$.  This is easily seen to have
finite-dimensional cohomology, and hence so does the extension
$(\Lambda a \otimes \mathcal{M}  \otimes \mathcal{N} , D)$.
Therefore, we have $\mathrm{rk}_0(\mathcal{M}\otimes \mathcal{N} )
\geq 1 > \mathrm{rk}_0(\mathcal{M} ) + \mathrm{rk}_0(\mathcal{N}
)$.
\end{example}

\begin{discussion}
Recall that $\mathrm{rk}(X)$ denotes the actual toral rank (rather
than its rational counterpart) for a finite complex $X$.  It is
easy to see that the inequality $\mathrm{rk}(X\times Y) \geq
\mathrm{rk}(X) + \mathrm{rk}(Y)$ holds in general---one takes the
obvious product action.  The minimal models of \exref{ex:Steve's
example} can be realised as the models of honest `geometrical'
spaces. $\mathcal{M} = \Lambda(x, y)$ is the minimal model of
$S^{12}$.  $\mathcal{N} = \Lambda(u_1, u_2, u_3, u_4, u_5, u_6, w,
v)$ is the minimal model of a sphere bundle over a product of odd
spheres.  We sketch this:  Take $S = S^{3}\times S^{3}\times
S^{3}\times S^{3}\times S^{5}\times S^{19}$, so that $S$ has
minimal model $\Lambda(u_1, u_2, u_3, u_4, u_5, u_6)$ with zero
differential.  $S$ has dimension $36$ and by pinching the
$35$-skeleton to a point, we obtain a non-trivial map $q \colon S
\to  S^{36}$.  Composing this with a suitable map $p\colon S^{36}
\to BSO(19)$ to the universal bundle $BSO(19)$, we obtain a
non-trivial map $S \to BSO(19)$.  Pulling back the universal real
oriented line bundle to one over $S$, via this map, we obtain a
real oriented line bundle $\mathbb{R}^{19} \to E \to S$.  Finally,
taking the unit sphere bundle gives a sphere bundle
\begin{displaymath}
S^{18} \to Y \to S.
\end{displaymath}
We assert that this sphere bundle $Y$ has minimal model
$\mathcal{N}$, as in \exref{ex:Steve's example}. Our assertion can
be justified by a consideration of the possible forms that the
minimal model of $Y$ can take, bearing in mind the construction of
$Y$.  Our main point here is simply to indicate how $\mathcal{N}$
corresponds to the minimal model of a reasonable space.  As in
\exref{ex:Steve's example}, $\mathrm{rk}_0(X) = \mathrm{rk}_0(Y) =
0$ and this is sufficient to conclude $\mathrm{rk}(X) =
\mathrm{rk}(Y) = 0$.  On the other hand, we have $\mathrm{rk}_0(X
\times Y) \geq 1$, but this is \emph{not} sufficient to conclude
$\mathrm{rk}(X \times Y) \geq 1$. We are left with the following
intriguing question: Does this space $X\times Y$ admit an
almost-free circle action? Note that by \cite[Prop.4.2]{Halperin},
there is a simply-connected, finite complex that admits a free
circle action and which has the rational homotopy type of $X
\times Y$.
\end{discussion}

The preceding examples and discussion give rise to a number of
interesting questions.  Generally, it would be useful to have
conditions under which the equality $\mathrm{rk}_0(X \times Y) =
\mathrm{rk}_0(X) + \mathrm{rk}_0(Y)$ either holds, or does not
hold. Various special cases are also of interest. For instance, we
can ask when does $\mathrm{rk}_0(X \times S^{2n+1}) =
\mathrm{rk}_0(X) + 1$?  We note that Halperin has an example in
which $\mathrm{rk}_0(X \times S^{2n+1}) > \mathrm{rk}_0(X) + 1$.
As another special case we could ask whether, for an $n$-fold
product of a space with itself, we have $\mathrm{rk}_0(X^n) =
n\,\mathrm{rk}_0(X)$? At present, we know of no example where
equality does not hold. Finally, we note that, as in the
discussion above, it is reasonable to ask all these questions in
the integral setting, too.

\subsection{The Gottlieb group} An interesting suggestion arising
from \thmref{thm:general 2-stage} is a connection between
$\mathrm{rk}_0(X)$ and the \emph{dimension of the Gottlieb group}.
We recall that the $n^{\text{th}}$ Gottlieb group is the subgroup
of $\pi_n(X)$ of elements $\a=[f] $  such that $(f,\id):\sph^n
\vee X \to X$ extends to a continuous map $\sph^n\times X \to X$.
In terms of a minimal model $(\Lam W,d)$, this corresponds
\cite[p.392]{F-H-T00} to the subspace
$$G_n(X)=\left\{
\begin{aligned}
 \beta\in (W^n)^* \mid \beta \text{ extends to a
derivation }&\theta\text{ of degree } -n \text{ of }(\Lambda W,d)\\
&\text{such that } d\theta - (-1)^n \theta d = 0
\end{aligned}
\right\}$$
In the two-stage case, with all generators of odd degree, the
minimal model can be written so that $\mathrm{dim}\,V =
\mathrm{dim}\,G_*(X)$. Thus we have the following observation:

\begin{corollary}
Suppose $X$ has minimal model that is two-stage with odd-degree
generators only.  Suppose that $G_*(X)$ denotes the rational
Gottlieb group of $X$. Then
\begin{displaymath}
\mathrm{dim}\,H(X) \geq 2^{\mathrm{dim}\,G_*(X)}.
\end{displaymath}
\end{corollary}

An examination of \exref{ex:two-stage odd} shows that we can have
$\mathrm{rk}_0(X) > \mathrm{dim}\,G_*(X)$. Next, we illustrate
that the dimension of the Gottlieb group can exceed the toral rank
by an arbitrary amount.

\begin{example}\label{ex:Gottlieb vs. rank}
For each $n \geq 1$, we describe an $(n+1)$-stage model
$\mathcal{M}_n$. Set $\mathcal{M}_n = \Lambda(x_1, x_2, y_1, z_1,
y_2, z_2, \ldots, y_{n-1}, z_{n-1}, y_n)$.  For degrees we have
$|x_i| = 3$ and $|z_i| = 3$ for each $i$.  The degrees of the
$y_i$ are chosen to be compatible with the differential, which is
as follows: $d(x_i) = 0$ and $d(z_i) = 0$ for each $i$, then
\begin{align*}
d(y_1) &= x_1x_2 \\ d(y_2) &= x_1x_2y_1z_1 \\ d(y_3) &=
x_1x_2y_1z_1y_2z_2 \\
  & \vdots \\
d(y_n) &= x_1x_2y_1z_1y_2z_2\cdots y_{n-1}z_{n-1}
\end{align*}
For each $n$, one can show that $\mathrm{dim}\,G_*(X) = n$, but
$\mathrm{rk}_0(X) = 1$. Further, a direct computation shows that
$\mathrm{dim}\,H(X) = (1/3)4^{n+1} + 2/3$.
\end{example}

These considerations suggest the following question:

\begin{question}\label{que:TRC Gottlieb}
Let $X$ be a finite complex with rational homotopy groups non-zero
in odd degrees only.  When is $\mathrm{dim}\,H(X) \geq
2^{\mathrm{dim}\,G_*(X)}$?
\end{question}

In spite of the apparent difficulty in establishing
\conjref{conj:TRC}, it would appear that the conjectured lower
bound of $2^{\mathrm{rk}_0(X)}$ underestimates the dimension of of
the cohomology, in many cases quite seriously.  This is supported
by examples such as \exref{ex:Gottlieb vs. rank}, \cite[Ex.
4.5]{Allday-Puppe85b} and analogous computations of the cohomology
of nilpotent Lie algebras, such as \cite{A-C-Jessup}. It is
possible that the lower bound suggested in \queref{que:TRC
Gottlieb} might give a closer estimate in some cases.

\providecommand{\bysame}{\leavevmode\hbox
to3em{\hrulefill}\thinspace}

\end{document}